\documentclass{fic-l}

\newtheorem{theorem}{Theorem}[section]
\newtheorem{lemma}[theorem]{Lemma}
\newtheorem{proposition}[theorem]{Proposition}

\theoremstyle{definition}

\theoremstyle{remark}
\newtheorem{remark}[theorem]{Remark}

\numberwithin{equation}{section}

\makeatletter
\@addtoreset{equation}{section}
\makeatother
\renewcommand{\theequation}{\thesection.\arabic{equation}}

\def\RO{\mathrm{O}}
\def\SO{\mathrm{SO}}
\def\apeir{\mathop{\rm apeir}\nolimits}
\def\rank{\mathop{\rm rank}}
\newcommand{\hole}{{\mkern2mu|\mkern2mu}}
\def\mix{\mathbin{\Diamond}}

\newcommand{\ve}{\varepsilon}

\newcommand{\vt}{\vartheta}
\newcommand{\sg}{\sigma}
\newcommand{\vp}{\varphi}
\newcommand{\Ga}{\mathnormal{\Gamma}}
\newcommand{\D}{\mathnormal{\Delta}}
\newcommand{\Ph}{\mathnormal{\Phi}}
\newcommand{\Ps}{\mathnormal{\Psi}}
\newcommand{\CC}{\mathcal{C}}
\newcommand{\CG}{\mathcal{G}}
\newcommand{\CP}{\mathcal{P}}
\newcommand{\CQ}{\mathcal{Q}}
\newcommand{\CR}{\mathcal{R}}
\newcommand{\BC}{\mathbb{C}}
\newcommand{\BE}{\mathbb{E}}
\newcommand{\BH}{\mathbb{H}}
\newcommand{\BR}{\mathbb{R}}
\newcommand{\BS}{\mathbb{S}}
\newcommand{\SJ}{\mathsf{J}}
\newcommand{\SK}{\mathsf{K}}
\newcommand{\SN}{\mathsf{N}}
\newcommand{\seg}{\{\mkern4mu\}}
\newcommand{\scl}[1]{\langle\mkern1mu{#1}\mkern1mu\rangle}
\newcommand{\fh}{{\textstyle \frac{5}{2}}}
\newcommand{\half}{{\textstyle \frac{1}{2}}}
\newcommand{\tth}{{\textstyle \frac{10}{3}}}

\newcommand{\ol}[1]{\overline{#1}}
\newcommand{\bpf}{\begin{proof}}
\newcommand{\epf}{\end{proof}}

\newcommand{\beq}{\begin{equation}}
\newcommand{\beql}[1]{\begin{equation} \label{#1}}
\newcommand{\eeq}{\end{equation}}
\newcommand{\bqy}{\begin{eqnarray*}}
\newcommand{\eqy}{\end{eqnarray*}}
\newcommand{\bry}{\begin{array}}
\newcommand{\ery}{\end{array}}
\newcommand{\blem}{\begin{lemma}}
\newcommand{\bleml}[1]{\begin{lemma} \label{#1}}
\newcommand{\elem}{\end{lemma}}
\newcommand{\bprop}{\begin{proposition}}
\newcommand{\bpropl}[1]{\begin{proposition} \label{#1}}
\newcommand{\eprop}{\end{proposition}}
\newcommand{\bthm}{\begin{theorem}}
\newcommand{\bthml}[1]{\begin{theorem} \label{#1}}
\newcommand{\ethm}{\end{theorem}}
\newcommand{\brem}{\begin{remark}}
\newcommand{\breml}[1]{\begin{remark} \label{#1}}
\newcommand{\erem}{\end{remark}}
\newcommand{\bitem}{\begin{itemize}} 
\newcommand{\eitem}{\end{itemize}}
\newcommand{\bpc}{\begin{picture}}
\newcommand{\epc}{\end{picture}}

\newcommand{\sectref}[1]{Section~\ref{#1}}

\newcommand{\propref}[1]{Proposition~\ref{#1}}
\newcommand{\thmref}[1]{Theorem~\ref{#1}}
\newcommand{\remref}[1]{Remark~\ref{#1}}

\newcommand{\Resq}[1]{\put(#1){\bpc(48,48)
\thicklines
\multiput(0,0)(0,48){2}{\line(1,0){48}}
\multiput(0,0)(48,0){2}{\line(0,1){48}}
\multiput(0,0)(48,0){2}{\circle*{5}}
\multiput(0,48)(48,0){2}{\circle*{5}}
\epc}}

\newcommand{\Horone}[1]{\put(#1){\bpc(48,0)
\thicklines
\put(0,0){\line(1,0){48}}
\multiput(0,0)(48,0){2}{\circle*{5}}
\epc}}

\newcommand{\Verone}[1]{\put(#1){\bpc(0,48)
\thicklines
\put(0,0){\line(0,1){48}}
\multiput(0,0)(0,48){2}{\circle*{5}}
\epc}}

\newcommand{\Rtri}[1]{\put(#1){\bpc(40,48)
\thicklines
\put(0,0){\line(5,3){40}}
\put(0,48){\line(5,-3){40}}
\put(0,0){\line(0,1){48}}
\put(40,24){\circle*{5}}
\multiput(0,0)(0,48){2}{\circle*{5}}
\epc}}

\newcommand{\REsq}[1]{\put(#1){\bpc(72,72)
\thicklines
\multiput(0,0)(0,72){2}{\line(1,0){72}}
\multiput(0,0)(72,0){2}{\line(0,1){72}}
\multiput(0,0)(72,0){2}{\circle*{6}}
\multiput(0,72)(72,0){2}{\circle*{6}}
\epc}}

\begin{document}

\title{Regular and chiral polytopes in low dimensions}
\author{Peter McMullen}
\address{University College London\\ Gower Street\\ 
London WC1E 6BT, England}
\email{p.mcmullen@ucl.ac.uk}
\author{Egon Schulte}
\address{Northeastern University\\ 
Boston, MA 02115, USA}
\email{schulte@neu.edu}
\thanks{The second author was supported in part by NSA-grant H98230-04-1-0116.}
\date{November 24, 2004}

\begin{abstract}
There are two main thrusts in the theory of regular and chiral polytopes:  the abstract, purely combinatorial aspect, and the geometric one of realizations.  This brief survey concentrates on the latter.  The dimension of a faithful realization of a finite abstract regular polytope in some euclidean space is no smaller than its rank, while that of a chiral polytope must strictly exceed the rank.  There are similar restrictions on the dimensions of realizations of regular and chiral apeirotopes.  From the viewpoint of realizations in a fixed dimension, the problems are now completely solved in up to three dimensions, while considerable progress has been made on the classification in four dimensions, the finite regular case again having been solved.  This article reports on what has been done already, and what might be expected in the near future.
\end{abstract}

\maketitle

\section{Introduction}
\label{intro}

Donald Coxeter's work on regular polytopes and groups of reflexions is often viewed as his most important contribution. At its heart lies a dialogue between geometry and algebra which was so characteristic for his mathematics (see, for example,  \cite{c_rp,c_rcp,cm}).  This  paper is yet more evidence for his lasting influence on generations of geometers.

In \cite{ms_rpo} (see also \cite[Sections~7E, 7F]{arp}), we classified completely all the faithfully realized regular polytopes and discrete regular apeirotopes in dimensions up to three.  Further, in \cite{m_rpfr}, the first author classified the regular polytopes and apeirotopes of maximal rank in each higher dimension, and showed that chiral polytopes could not have full rank.  Last, in \cite{s_cp1,s_cp2}, the second author has found all the chiral apeirohedra in three dimensions.  

The present paper surveys the developments on realizations of regular or chiral polytopes, which have occurred since the publication of our book~\cite{arp}. There are two quite different ways to approach realizations.  The first, for which a fairly complete theory exists (at least, in the finite case), asks for a description of the space of all realizations (a kind of ``moduli space") of a given abstract regular polytope or apeirotope, with rank playing only a minor r\^ole (see \cite[Sections~5B, 5C]{arp} for further details).  The second, about which much less is known in general terms, asks for a classification of the realizations of all these polytopes and apeirotopes in a euclidean space of given dimension (in this case, it is usual to impose conditions such as faithfulness and discreteness).  This problem is solved in three dimensions.  The finite regular polyhedra have long been known;  adding to the Petrie-Coxeter apeirohedra of \cite{c_rsp}, Gr\"unbaum \cite{g_on} found all but one of the remaining regular apeirohedra, while Dress \cite{d1,d2} found the missing example, and proved that the classification was then complete.  We refer the reader to \cite{ms_rpo} for a quick method of arriving at the full characterization, including a discussion of the geometry of the regular apeirohedra and presentations of their symmetry groups, as well as for the enumeration of the regular $4$-apeirotopes in three dimensions.

In four dimensions, the currently open problems are those of classifying the finite regular polyhedra, and the regular apeirohedra and $4$-apeirotopes;  \cite{m_rpfr} solves the problems of the regular $4$-polytopes and $5$-apeirotopes.  The paper \cite{m_fdrp} in preparation actually settles the first of these problems (the polytopes with planar faces were classified in \cite{abm,bracho});  however, the other two, together with the corresponding classification problems for chiral polytopes, are still open, although some progress has been made on them.

\section{Regular and chiral polytopes}
\label{regchirpol}

For the general background on abstract regular polytopes, we refer the reader to the recently published monograph \cite{arp};  for the most part, we shall not cite original papers directly.  In this paper, we largely concentrate on the geometric aspects of the theory, that is, on realizations of regular polytopes.  

However, we begin with the more combinatorial picture.  An \emph{abstract polytope} of \emph{rank} $n$, or simply an (\emph{abstract}) \emph{$n$-polytope}, is a partially ordered set $\CP$ with a strictly monotone rank function, taking values in $\{-1,0,\ldots,n\}$.  The elements of rank $j$ are the \emph{$j$-faces} of $\CP$, or \emph{vertices}, \emph{edges} and \emph{facets} of $\CP$ if $j = 0$, $1$ or $n-1$, respectively. The maximal chains are the \emph{flags} of $\CP$ and contain exactly $n + 2$ faces, including a unique minimal face and a unique maximal face (usually omitted from the notation). Two flags are called \emph{adjacent} if they differ by one element; then $\CP$ is \emph{strongly flag-connected}, meaning that, if $\Ph$ and $\Ps$ are two flags, then they can be joined by a sequence of successively adjacent flags $\Ph = \Ph_0,\Ph_1,\ldots,\Ph_k = \Ps$, each of which contains $\Ph \cap \Ps$. Finally, if $F$ and $G$ are a $(j-1)$-face and a $(j+1)$-face with $F < G$, then there are exactly \emph{two} $j$-faces $H$ such that $F < H < G$.  An $n$-polytope $\CP$ is then called \emph{regular} if its combinatorial automorphism group $\Ga(\CP)$ (preserving the partial ordering) is (simply) transitive on its flags; in this case, if $\Ph$ is a (fixed) \emph{base flag} and, for $j = 0,\ldots,n-1$, $\rho_j$ is the automorphism which maps $\Ph$ to the adjacent flag $\Ph^{j}$ with a different $j$-face, then $\Ga(\CP)$ is generated by $\rho_{0},\ldots,\rho_{n-1}$. 

We can adopt (see \cite[Theorem~2E11]{arp}) the viewpoint that an abstract regular polytope is to be identified with its group.  The latter is precisely what is called a \emph{string C-group}; here, the ``C" stands for ``Coxeter", though not every C-group is a Coxeter group.  A string C-group $\Ga$ is a group generated by $n$ involutions $\rho_j$ (the \emph{distinguished generators}) with $j \in \SN := \{0,\ldots,n-1\}$, such that $\rho_j$ and $\rho_k$ commute if $0 \leq j \leq k - 2 \leq n-3$, and
\beql{intprop}
\scl{\rho_i \mid i \in \SJ} \cap \scl{\rho_i \mid i \in \SK} = \scl{\rho_i \mid i \in \SJ \cap \SK}
\eeq
for each $\SJ,\SK \subseteq \SN$;  the last is the \emph{intersection property}.  Each string C-group $\Ga$ then determines (uniquely) a regular $n$-polytope $\CP$ with $\Ga = \Ga(\CP)$. The \emph{$j$-faces} of $\CP$ are the right cosets $\Ga_j\sg$ of the \emph{distinguished subgroup}
\[  \Ga_j := \scl{\rho_i \mid i \neq j}  \]
for each $j \in \SN$, and two faces are incident just when they intersect (as cosets).  In fact, incidence actually induces an order relation:
\[  \Ga_j\sg \leq \Ga_k\tau \iff \Ga_j\sg \cap \Ga_k\tau \neq \emptyset \mbox{ and } j \leq k.  \]
Formally, we also adjoin two copies of $\Ga$ itself, as the (unique) $(-1)$- and $n$-faces of $\CP$. The maximal chains (with respect to this ordering) are the \emph{flags} of $\CP$;  the group $\Ga$ is then simply transitive on the flags of $\CP$.  In particular,  for $j = 0,\ldots,n-1$ the distinguished generator $\rho_j$ of $\Ga$ takes the \emph{base flag} $\Ph := \{\Ga_{-1},\Ga_0,\Ga_1,\ldots,\Ga_{n-1},\Ga_n\}$ into the adjacent flag $\Ph^j$ which differs from it in $\Ga_j$. Note that the distinguished subgroups $\Ga_{n-1} = \scl{\rho_0,\ldots,\rho_{n-2}}$ and $\Ga_0 = \scl{\rho_1,\ldots,\rho_{n-1}}$ are themselves string C-groups;  the corresponding polytopes are the \emph{facet} and \emph{vertex-figure} of $\CP$, respectively (the latter consists of the faces of $\CP$ with vertex $\Ga_{0}$). As we said earlier, \cite[Theorem~2E11]{arp} shows that this description of a regular polytope $\CP$ in terms of (its C-group) $\Ga(\CP)$ and the previous one in terms of the face poset are equivalent.

The distinguished generators of $\Ga = \Ga(\CP)$ satisfy relations
\beq
\label{cgp}
(\rho_{i}\rho_{j})^{p_{ij}} = \ve \quad  (i,j = 0,\ldots,n-1),
\eeq
with $p_{ii} = 1$, $p_{ij} = p_{ji} \geq 2$ if $i \neq j$, and $p_{ij} = 2$ if $|i-j| \geq 2$ (hence the term ``string" C-group). The numbers $p_{j} := p_{j-1,j}$ ($j = 1,\ldots,n-1$) determine the \emph{Schl\"afli type} $\{p_{1},\ldots,p_{n-1}\}$ of $\CP$.  To avoid cases which, in our context, turn out to be trivial, we always assume that adjacent generators $\rho_{j-1}$ and $\rho_j$ of $\Ga$ do not commute (this is justified in \sectref{real}); in other words, $p_j > 2$ (possibly, $p_j = \infty$).   If the polytope is determined just by the $p_j$, then we have the \emph{universal} regular polytope (of that Schl\"afli type), for which we use the same symbol $\{p_1,\ldots,p_{n-1}\}$ (but without qualification);  we write $[p_1,\ldots,p_{n-1}]$ for the corresponding \emph{Coxeter} group.  Generally, however, the group $\Ga$ will satisfy additional relations as well, for some of which we introduce special notation later.

The underlying face-set of a polytope $\CP$  can be finite or infinite. An infinite $n$-polytope is also called an (\emph{abstract}) \emph{$n$-apeirotope}; when $n = 2$, we also refer to it as an \emph{apeirogon}, and when $n = 3$ as an \emph{apeirohedron}.  

A central question in the abstract theory is that of the amalgamation of polytopes of lower rank. If a regular $(n+1)$-polytope has facets (of type) the $n$-polytope $\CP$ and vertex-figures the $n$-polytope $\CQ$, then the facets of $\CQ$ must be isomorphic to the vertex-figures of $\CP$.  Conversely, if $\CP$ and $\CQ$ satisfy this latter criterion, then we write $\scl{\CP,\CQ}$ for the class of all regular $(n+1)$-polytopes with facet $\CP$ and vertex-figure $\CQ$.  The question has two parts.  First, is $\scl{\CP,\CQ} \neq \emptyset$;  in other words, does there exist any such regular $(n+1)$-polytope at all?  If so, then there is a \emph{universal} member $\{\CP,\CQ\}$ in the family $\scl{\CP,\CQ}$, of which every other one is a quotient (in the sense that its group is an appropriate quotient).  Second, given that it exists, we ask what $\{\CP,\CQ\}$ is.  (See \cite[Section~4B]{arp} for further details.)  In the present context, we often pose this question in the form:  is a given regular polytope, whose facet and vertex-figure are known, actually universal of its kind?

There are several general techniques for constructing new regular polytopes from old ones. In particular, two different regular polytopes may be related by what is called a \emph{mixing operation};  the distinguished generators of the second group are certain products of those of the first (see \cite[Chapter~7]{arp}).  Apart from the \emph{duality} operation $\delta$, which just reverses the order of the distinguished generators (and the order relation on the faces), there are two others we mention here; one further operation (for chiral polyhedra) will occur in Section~\ref{chirpol3}.  Let $\Ga = \scl{\rho_i \mid i \in \SN}$ be a string C-group, let $j \neq k$, and consider the operation
\[  (\rho_0,\ldots,\rho_{n-1}) \mapsto (\rho_0,\ldots,\rho_{j-1},\rho_j\rho_k,\rho_{j+1},\ldots,\rho_{n-1}) =: (\sg_0,\ldots,\sg_{n-1}).  \]
Since adjacent generators of $\Ga$ do not commute, we easily see that the group $\D := \scl{\sg_0,\ldots,\sg_{n-1}}$ cannot possibly be a string C-group unless $(j,k) = (2,0)$ or $(n-3,n-1)$.  The former will rule itself out later for geometric reasons (see Section~\ref{real});  the latter, namely,
\beql{petrie}
\pi\colon\ (\rho_0,\ldots,\rho_{n-1}) \mapsto (\rho_0,\ldots,\rho_{n-4},\rho_{n-3}\rho_{n-1},\rho_{n-2},\rho_{n-1}) =: (\sg_0,\ldots,\sg_{n-1}),
\eeq
which we denote by $\Ga \mapsto \Ga^\pi$, is called the \emph{Petrie operation}, since it generalizes the operation with the same name when $n = 3$.  Even when $n = 3$, the Petrie operation $\pi$ does not always yield a C-group (though such cases are rather exceptional), but, for higher rank, each application has to be checked directly. However, if in fact $\Ga^\pi$ is a C-group, then we write $\CP \mapsto \CP^\pi$ to indicate the effect of the operation on the corresponding polytope $\CP$;  the new polytope $\CP^\pi$ is called the \emph{Petrial} of $\CP$.  One general case (see \cite{m_rpfr}) can be settled easily.
\bpropl{nonpetrie}
If $\Ga = \scl{\rho_0,\ldots,\rho_{n-1}}$ is a string C-group with $n \geq 4$ for which $p_{n-3}$ is odd, then the Petrial $\Ga^{\pi}$ is not a C-group.
\eprop

Mixing operations are particularly powerful when applied to regular polyhedra or apeirohedra $\CP$. For example, the Petrial $\CP^{\pi}$ can be obtained from $\CP$ by replacing the $2$-faces by the \emph{Petrie polygons} of $\CP$ (while keeping the vertices and edges); the geometric picture of a Petrie polygon here is one which shares two successive edges of each $2$-face which it meets, but not a third. An important class of regular polyhedra or apeirohedra consists of those which are completely determined by their Schl\"afli type and the length of their Petrie polygons.   We write $\{p,q\}_r$ for the polyhedron (possibly infinite) of Schl\"afli type $\{p,q\}$, whose Petrie polygons of length $r$ determine it.  Its group is the Coxeter group $\scl{\rho_0,\rho_1,\rho_2} = [p,q]$, with the imposition of the single extra relation
\beql{petriepol}
(\rho_0\rho_1\rho_2)^r = \ve.
\eeq
We note that, if it is a genuine polyhedron, then the Petrial of $\{p,q\}_r$ is $\{r,q\}_p$. 

In the context of polyhedra, another operation is also of great importance. The (\emph{second}) 
\emph{facetting operation} $\vp_{2}$ is given by
\beql{facett}
\vp_{2}\colon\  (\rho_{0},\rho_{1},\rho_{2})  \mapsto (\rho_{0},\rho_{1}\rho_{2}\rho_{1},\rho_{2}), 
\eeq
and replaces the $2$-faces of a polyhedron $\CP$ by the \emph{holes} (while keeping the vertices and edges);  a hole of $\CP$ is an edge-circuit which exits from the \emph{second} edge (in some local orientation) emanating from a vertex from the edge by which it entered.  The designation of a (possibly infinite) regular polyhedron of Schl\"afli type $\{p,q\}$, which is determined by its holes of length $h$, is 
$\{p,q \hole h\}$.  The corresponding relation to be imposed on the Coxeter group $\scl{\rho_0,\rho_1,\rho_2} = [p,q]$ is 
\beql{holepol}
(\rho_0\rho_1\rho_2\rho_1)^h = \ve.
\eeq
Various examples of such polyhedra occur later;  for now, let us observe that the three Petrie-Coxeter apeirohedra are, as abstract regular polyhedra, $\{4,6 \hole 4\}$, $\{6,4 \hole 4\}$ and $\{6,6 \hole 3\}$. 

In \cite[Section 7A]{arp} we also introduced the notion of a mix of two regular polytopes (or corresponding C-groups).  The following abstract construction is a special case of this mix and occurs when one polytope is $1$-dimensional, that is, a segment.  Again, suppose that $\Ga = \scl{\rho_i \mid i \in \SN}$ is a string C-group.  Let $\tau$ be an involution which commutes with all $\rho_j$, and consider the operation
\beql{mixseg}
(\rho_0,\ldots,\rho_{n-1},\tau) \mapsto (\rho_0\tau,\rho_1\ldots,\rho_{n-1}) =: (\sg_0,\ldots,\sg_{n-1}).
\eeq
This is called \emph{mixing with a segment}, because $\tau$ can be regarded as the generating involution of the group of the segment $\seg$ (see \sectref{real} for the notation).  We have (see \cite[Theorem~7A8]{arp})
\bthml{mixsegpol}
Mixing a string C-group $\Ga$ with the group of a segment always yields another C-group.  This is isomorphic to $\Ga$ if all edge-circuits in the associated regular polytope $\CP$ have even length;  otherwise, it is isomorphic to the direct product $\Ga \times \CC_2$ of $\Ga$ with a cyclic group $\CC_2$ of order $2$.
\ethm

The resulting regular polytope (which we again say is obtained from $\CP$ by mixing with a segment) is denoted by $\CP \mix \seg$.  This has twice as many vertices as $\CP$ precisely when some edge-circuit of $\CP$ has odd length.

We also require another basic technique for constructing regular polytopes from certain groups by what are called \emph{twisting operations} (see \cite[Chapter~8]{arp}).  In this, a given group (usually itself a C-group) is augmented by means of one or more group automorphisms.  This technique has been extremely successful in various classification problems for regular polytopes.  In the present context, it assumes great importance in the enumeration of the regular polyhedra in $\BE^4$;  see \sectref{regpol4} below.

Roughly speaking, chiral polytopes have half as many possible automorphisms as have regular polytopes.  More technically, the $n$-polytope $\CP$ is \emph{chiral} if it has two orbits of flags under its group $\Ga(\CP)$, with adjacent flags in different orbits.  A chiral $n$-polytope $\CP$ is then identified with a group of the form $\Ga = \scl{\sg_1,\ldots,\sg_{n-1}}$, on which there are relations
\beql{chirpolrel}
\begin{cases}
\sg_j^{p_j} = \ve, & \text{$j = 1,\ldots,n-1$,} \cr
(\sg_j \sg_{j+1} \cdots \sg_k)^2 = \ve, & \text{$1 \leq j < k \leq n-1$.}
\end{cases}
\eeq
We again refer to $\{p_1,\ldots,p_{n-1}\}$ as the \emph{Schl\"afli type} of $\CP$.

The relationship between the group and the corresponding (abstract) polytope is a little less obvious than is the case for regular polytopes (see \cite{sw_c} for more details).  The distinguished generator $\sg_j$ permutes the $(j-1)$- and $j$-faces cyclically in the appropriate section of the base flag $\Ph = \{F_0,F_1,\ldots,F_{n-1}\}$;  if $F_j'$ replaces $F_j$ in the adjacent flag $\Ph^j$, then $F_{j-1}'\sg_j = F_{j-1}$ and $F_j\sg_j = F_j'$.  The vertices of $\CP$ are (identified with) the right cosets of the subgroup $\Ga_0 := \scl{\sg_2,\ldots,\sg_{n-1}}$, with $F_0 = \Ga_0$ itself the base vertex.  The involutory element $\tau := \sg_1\sg_2$ interchanges the two vertices of the base edge, taking $\Ph$ into $(\Ph^{0})^{2} = (\Ph^{2})^{0}$;  it is often useful to replace $\sg_1$ as a generator by $\tau$ (compare \cite{s_cp1}).

In a chiral polytope, adjacent flags are not equivalent under the group. If $\Ph$ is replaced by an adjacent flag, $\Ph^{0}$ (say), then the respective generators are $\sg_{1}^{-1}, \sg_{1}^{2}\sg_{2}, \sg_{3}, \ldots,\sg_{n-1}$. Thus a chiral polytope occurs in two ({\em combinatorially}) {\em enantiomorphic forms\/}, each specified by the choice of an orbit of base flags ($\Ph$ or $\Ph^{0}$), or, equivalently, a conjugacy class of sets of generators (represented by $\sg_{1},\ldots,\sg_{n-1}$ or $\sg_{1}^{-1}, \sg_{1}^{2}\sg_{2}, \sg_{3}, \ldots, \sg_{n-1}$, respectively).  For a regular polytope, these two enantiomorphic forms can be identified (under the generator $\rho_{0}$ of $\Ga$).

\section{Realizations}
\label{real}

There are many candidates for spaces in which regular polytopes $\CP$ might be realized geometrically. The usual  (and generally most useful) context of realizations is of those in euclidean spaces, because it is in these that we obtain the richest structure.  However, initially at least, it is appropriate for us to broaden the definition.  Thus, for the time being, $E$ is a $k$-dimensional spherical space $\BS^k$, euclidean space $\BE^k$ or hyperbolic space $\BH^k$, for some $k$. If $\CP$ is a finite polytope, then $E$ will be spherical;  if $\CP$ is an apeirotope, then, since we are generally interested only in discrete realizations, $E$ will be euclidean or hyperbolic. 

We begin with a brief review of some definitions (see \cite[Chapter~5]{arp} for the general background here). Let $\CP$ be an abstract regular polytope (or apeirotope -- for the moment, we use the generic term, not distinguishing between the finite and infinite cases), and let $\Ga:=\Ga(\CP)$. For a \emph{faithful realization} of $\CP$ we have two ingredients. First, we need a suitable space $E$ which admits a group $\CG$ of isometries isomorphic to $\Ga$; this is the \emph{symmetry group} of the realization of $\CP$. It is convenient to identify the \emph{reflexion} $R_j$ in $\CG$ corresponding to the involution $\rho_j$ in $\Ga$ with its \emph{mirror} 
\[  \{x \in E \mid xR_j = x\}  \]
of fixed points;  we thus use the same symbol $E$ for the ambient space to denote the identity mapping. The intersection
\[  W := R_1 \cap \cdots \cap R_{n-1}  \]
is called the \emph{Wythoff space} of the realization.  The realization of $\CP$ associated with
$\CG$ and its generators $R_{j}$ then arises from some choice of \emph{initial vertex} $v \in W$.  The vertex-set of the realization is $V := v\CG$, the orbit of $v$ under $\CG$, and we always assume that $E$ is spanned by $V$ (as a subspace of the appropriate kind), so that $E$ is thought of as the \emph{ambient space} of the realization, namely, the space (of one of the three kinds) of smallest dimension which contains it. 

Note that, if $\CG$ were to be such that $R_j = E$, the identity mapping, then $R_k = E$ for all $k > j$ as well and the realization would not be faithful.  In particular, this will happen if $p_j = 2$, which is why we excluded this possibility in \sectref{regchirpol}.

The induced geometric structure, the actual \emph{realization} $P$ of $\CP$, is defined as follows. Write $F_0 := v$, and, for $j\geq 1$, let
\[  F_j := F_{j-1}\scl{R_0,\ldots,R_{j-1}};  \]
these are the basic faces. Then the $j$-faces of the realization are the $F_jG$ with $G\in\CG$, with the order relation given by iterated membership. Thus \emph{edges} are composed of the two vertices which belong to them (we also think of an edge as the line-segment between its vertices -- there will be no ambiguity, even in the spherical case, because antipodal points of the sphere will never determine an edge), $2$-faces of the edges which belong to them, and so on up to the \emph{ridges} or $(n-2)$-faces and \emph{facets} or $(n-1)$-faces. We sometimes refer to the realization $P$ as a \emph{geometric polytope}. Its \emph{dimension} is defined by $\dim P := \dim E$, and its vertex-set is denoted by $V(P) := V$. Finally, for the realization to be \emph{faithful\/}, we demand that, for each $j = 1,\ldots,n-1$, a $j$-face be uniquely determined by the $(j-1)$-faces which belong to it. Recall here our initial assumption that $\CG$ and $\Ga$ be isomorphic, so for a faithful realization we then have natural bijections between the sets of $j$-faces of $\CP$ and $P$ for each $j$. Some regular polytopes do not admit faithful realizations, because this latter condition implies a corresponding purely combinatorial condition on $\CP$.

A realization of an abstract regular $n$-polytope $\CP$ determines a realization of each of its faces or co-faces (iterated vertex-figures). In particular, $F_{n-1}$ (and its induced structure, with the same initial vertex $v$) gives a realization of the facet of $\CP$;  its symmetry group is the image $\CG_{n-1}$ of $\Ga_{n-1}$.  If we write $w$ for the mid-point of the edge between $v$ and $vR_0$, then $w$ is the initial vertex of a realization of the vertex-figure of $\CP$, with symmetry group the image $\CG_0$ of $\Ga_0$.  (This suffices for our purposes. However, in the hyperbolic case of a polytope with vertices on the absolute, then the initial vertex $w$ is well-defined as the intersection of the mirror $R_0$ with the line between $v$ and $vR_0$ -- in any event, $w$ will always lie in this intersection.)  Faithfulness is hereditary; that is, if the original realization of $\CP$ is faithful, then the realizations of the facet and vertex-figure of $\CP$ are also faithful.  In a similar way, $\scl{R_0,\ldots,R_{j-1}}$ is the symmetry group of the basic $j$-face $F_j$ of $P$, while $\scl{R_{j+1},\ldots,R_{n-1}}$ is that of the basic \emph{co-$j$-face} $P/F_j$, which is the $(j+1)$-fold iterated vertex-figure.  Thus the vertex-figure itself is $P/F_0$.  Even more generally, $\scl{R_{j+1},\ldots,R_{k-1}}$ is the symmetry group of the \emph{section} $F_k/F_j$ (for $j \leq k-2$), the $(j+1)$-fold iterated vertex-figure of the basic $k$-face $F_{k}$. 

We often find it more convenient to use $vR_0$ rather than $w$ as the initial vertex of the
vertex-figure;  for most purposes, this makes little difference, since the combinatorics are not altered.

For regular polytopes of rank at most $2$ we have the following spherical or euclidean realizations.  In $\BE^0$ we just have the point (realizing the $0$-polytope), the finite regular $1$-polytopes are segments $\seg$, which are naturally realized in the $0$-sphere $\BS^0$, while the regular apeirogon $\{\infty\}$ is naturally realized discretely in $\BE^1 = \BR$.  In the unit circle $\BS^1$, there is an infinite family of (finite) regular polygons.  Their mirrors $R_0$ and $R_1$ are lines through its centre at a \emph{rational} angle $\pi/p$, meaning that $p > 2$ is a rational number (always in its lowest terms);  the resulting regular polygon is denoted $\{p\}$.  In addition, $\{\infty\}$ has non-discrete faithful realizations in $\BS^1$.  As we mentioned before, we shall not address here the question of finding all possible realizations of a given abstract regular polytope; a fairly complete theory has been described in \cite[Sections~5B, 5C]{arp}.  Suffice it to remark that the realization space has been determined for several interesting classes of polytopes; see, for example, \cite{mowe}.

There are important restrictions on faithful realizations;  we refer to \cite[Sections~5B, 5C]{arp} for proofs.
\bthml{rankdimpol}
Let $P$ be a faithful realization of an abstract regular polytope $\CP$, whose ambient space $E$ is a spherical, euclidean or hyperbolic space.  Then $\dim P \geq \rank\CP - 1$.
\ethm

\bthml{dimmirror}
Let $P$ be a faithful realization of an abstract regular $n$-polytope in $E$, with group $\CG = \scl{R_0,\ldots,R_{n-1}}$.  Then $\dim R_j \geq j$ for $j = 0,\ldots,n-2$, and $\dim R_{n-1} \geq n-2$.
\ethm

In both theorems, if the polytope is finite, so that the ambient space is spherical, then, regarded as euclidean realizations, each of the dimensions must be increased by $1$. 

If we have (not necessarily faithful) realizations of the abstract regular polytope (or apeirotope) $\CP$ in two euclidean spaces, say $P$ with mirrors $S_0,\ldots,S_{n-1}$ in $L$ and $Q$ with mirrors $T_0,\ldots,T_{n-1}$ in $M$ (possibly some $S_j = L$ or $T_j = M$), then their \emph{blend} has mirrors $S_j \times T_j$ in $L \times M$ for $j = 0,\ldots,n-1$.  Indeed, if $v \in S_1 \cap \cdots \cap S_{n-1}$ and $w \in T_1 \cap \cdots \cap T_{n-1}$ are the initial vertices of the two realizations, then $(v,w)$ can be chosen as the initial vertex of the blend, which we then write $P \# Q$.  A realization which cannot be expressed as a blend in a non-trivial way is called \emph{pure}.

One main tool for classifying regular polytopes of a fixed rank $n$ in a fixed dimension is the \emph{dimension vector} $(\dim R_0,\dim R_1,\ldots,\dim R_{n-1})$ of the possible realizations;  the first step in any enumeration is to determine which dimension vectors can occur.

It is worth noting that, in general, duals of faithfully realizable regular polytopes are not necessarily faithfully realizable at all (Petrials are particular examples), let alone in the same space.

There is a similar realization theory for chiral polytopes. Indeed, let us call a realization $P$ of an abstract polytope $\CP$ \emph{chiral} if $P$ has two orbits of flags under its symmetry group $\CG(P)$, with adjacent flags lying in different orbits.  It is clear that the original polytope $\CP$ must be regular or chiral. Note that there exist (already in $\BE^3$) faithful realizations of polytopes with two flag orbits under $\CG(P)$ which are not chiral (see \cite{w_i2} for examples).  

It is helpful to remark that, if $\CP$ is a regular $n$-polytope with group $\Ga = \scl{\rho_0,\ldots,\rho_{n-1}}$, then its combinatorial rotation subgroup $\Ga^+(\CP)$ has generators
\[  \sg_j := \rho_{j-1}\rho_j, \qquad j = 1,\ldots,n-1.  \]
Thus a chiral realization of a polytope may be thought of as having only rotational symmetries.  Moreover, if the abstract polytope $\CP$ is at least chiral, in that its group $\Ga$ contains the automorphisms $\sg_1,\ldots,\sg_{n-1}$ in the definition of chirality, then $\CP$ is actually regular if we can adjoin any one of the involutions $\rho_j$ for $j = 0,\ldots,n-1$.  (We then have $\rho_i = \sg_{i+1}\rho_{i+1}$ for $i = 0,\ldots,j-1$, or $\rho_i = \rho_{i-1}\sg_i$ for $i = j+1,\ldots,n-1$.)

Chiral realizations are derived by a variant of Wythoff's construction, applied to a suitable representation $\CG = \scl{S_1,\ldots,S_{n-1}}$ of the underlying combinatorial group $\Ga := \scl{\sg_{1},\ldots,\sg_{n-1}}$; the latter is $\Ga(\CP)$ or $\Ga^{+}(\CP)$ according as the abstract polytope $\CP$ is chiral or regular. The Wythoff space now is the fixed set of the subgroup $\CG_0 := \scl{S_2,\ldots,S_{n-1}}$. We describe the $3$-dimensional case in more detail in Section~\ref{chirpol3}.

It is clear that an abstract regular polytope may have chiral realizations, though not necessarily faithful ones;  it is an interesting open question whether it could actually have faithful chiral realizations.  It is an elementary observation that a realized polygon with full rotational symmetry group is actually regular.  Similar arguments to those used in the proof of Theorem~\ref{rankdimpol} then yield
\bpropl{rankchirreal}
If $P$ is a faithful chiral realization of an abstract polytope, whose ambient space is a spherical, euclidean or hyperbolic space $E$, then $\dim P \geq \rank \CP - 1$.
\eprop

When the abstract polytope $\CP$ is finite, we usually assume that the centroid of the vertex-set $V$ of its (chiral or regular) realization $P$ is the origin $o$ of $E$, so that $\CG$ is an orthogonal group.  If $\CP$ is infinite, in which case we again call $P$ a (\emph{geometric}) \emph{apeirotope}, we will additionally demand of $P$ that it be discrete, so that the group $\CG$ acts discretely on the ambient space $E$.  Moreover, in order to avoid constant repetition of various fixed phrases subsequently, we adopt the conventions that, in the geometric context of realizations, \emph{regular polytope} will mean ``faithfully realized finite abstract regular polytope'', while \emph{regular apeirotope} will mean ``discrete faithfully realized abstract regular apeirotope''; we also adopt the corresponding terminology for chiral polytopes and chiral apeirotopes.

We end the section with two general remarks.  Let $S$ and $T$ be linear reflexions.  First, since $ST = (-S)(-T) = S^\perp T^\perp$ (thus identifying $-S$ with its mirror $S^\perp$, and so on), then $S \cap T$ and $S^\perp \cap T^\perp$ are both pointwise fixed by the product.  That is, the axis (fixed set) of $ST$ is
\beql{prodaxis}
(S \cap T) + (S^\perp \cap T^\perp)\;\, (\,= (S \cap T) + (S+T)^\perp\,).
\eeq
In particular, if $S$ and $T$ commute, then \eqref{prodaxis} is the mirror of their product $ST = TS$, which is again a reflexion.

Second, we have a general construction from \cite{m_rpfr}, of which special cases already occur in \cite{ms_rpo}.  Let $X$ be a point-set in a euclidean space $E$.  We call $X$ \emph{rational} if the points of $X$ can be chosen to have rational coordinates with respect to some (linear or affine) coordinate system in $E$.  The following remark is obvious.
\bleml{apeirset}
Let $E$ be a euclidean space, and let $X$ be a finite point-set in $E$.  Let $\CR(X)$ be the group generated by the point-reflexions (inversions) in the points of $X$.  Then $\CR(X)$ is discrete if and only if $X$ is rational.
\elem

If $P$ is a regular polytope with ambient space $E$, then we similarly call $P$ \emph{rational} if its vertex-set is rational.  We have the following.
\bthml{apeirpol}
Let $P$ be a rational regular $n$-polytope in the euclidean space $E$, with symmetry group $\CG_0 = \scl{R_1,\ldots,R_n}$ and initial vertex $w$, and suppose that $v \in R_1 \cap \cdots \cap R_n$.  Let $R_0 = \{w\}$ be the point-reflexion in the point $w$.  Then $\CG := \scl{R_0,\ldots,R_n}$ is the group of a discrete regular $(n+1)$-apeirotope $\apeir P$, with $2$-faces apeirogons, and vertex-figure $P$ at the initial vertex $v$.
\ethm

We call $\apeir P$ the \emph{free abelian apeirotope on} $P$, or \emph{with vertex-figure} $P$, and \emph{base vertex} $v$.  When we apply this construction, it will usually be the case that $P$ itself is finite and full-dimensional in $E$, so that $v$ is the centre of $P$.

\section{Regular polytopes of full rank}
\label{fullrank}

If $P$ is a realization of a regular polytope $\CP$ for which equality holds in \thmref{rankdimpol}, then we say that $P$ is \emph{of full rank}.  The emphasis is placed this way round, because our aim (as explained in \sectref{intro}) is to classify regular (and chiral) polytopes by dimension.  In this case, we can go further than \thmref{dimmirror}, and place further restrictions on the dimensions of the mirrors of the generating reflexions of the realizations.  We refer to \cite{m_rpfr} for a proof.
\bthml{fullrankmir}
Let $P$ be a faithful realization of full rank of a regular $n$-polytope $\CP$ in the ambient space $E$, with symmetry group $\CG = \scl{R_0,\ldots,R_{n-1}}$.  Then $\dim R_j = j$ or $n-2$ for $j = 0,\ldots,n-3$, and $\dim R_{n-2} = \dim R_{n-1} = n-2$.
\ethm

For finite polytopes, we now find it convenient to revert to the former definition of realization in euclidean spaces.  In other words, henceforth we regard a sphere which carries the vertices of a realization $P$ of a finite regular polytope as sitting in the euclidean space of one larger dimension with centre the origin $o$.  The mirrors $R_j$ of its euclidean group $\CG$ are then thought of as linear subspaces, also of one larger dimension than before;  in particular, in the minimal case, $R_0$ is either a line or a hyperplane.  Finally, we shall use the more familiar $I$ for the identity (in a sense, $E$ is no longer quite appropriate).

\breml{reflgpquery}
If $R$ is a linear reflexion in a euclidean space $E$, then  $-R = (-I)R$, the product of $R$ with the central inversion $-I$, is the reflexion in the orthogonal complement $R^\perp$ of $R$.  Replacing a mirror by its orthogonal complement is often a useful tool in studying realizations. In particular, in the case of a faithful realization of full rank of a finite regular $n$-polytope with centre $o$, if the mirror $R_0$ is a line, then $-R_0$ is a hyperplane reflexion.  If we replace $R_0$ by $-R_0$, then at worst we have replaced the symmetry group $\CG$ by $\CG \times \CC_2$, with $\CC_2 = \{\pm I\}$;  in any event, we always have another finite group. Thus the mirror replacement often produces groups closely related to finite groups generated by hyperplane reflexions.
\erem

\remref{reflgpquery} enables us to introduce some important geometric operations on finite polytopes of full rank, which are the key to their enumeration.  For such polytopes $P$, since $o$ is the sole fixed point of the ambient space $E$ under the group $\CG$, it follows that 
\[ K_0 := R_0 \cap \cdots \cap R_{n-1} = \{o\}. \] 
Thus the central reflexion $-I$, identified with its mirror $\{o\}$, is $K_0$, so the mirror replacement  of \remref{reflgpquery} is $R_0 \mapsto R_0K_0$.  Moreover, it is extremely useful to have variant operations, which act on the co-$(j-1)$-face $P/F_{j-1}$ for some $j$ and also apply to apeirotopes when their co-$(j-1)$-faces are finite. With
\[  K_k := R_k \cap \cdots \cap R_{n-1} \quad (0 \leq k \leq n-1),  \]
we see that (the reflexion in) $K_k$ induces the central inversion on the affine hull of $P/F_{k-1}$;  recall here our assumption of full rank.  For $0 \leq j \leq k \leq n-1$, we then define the operation $\kappa_{jk}$ on $\CG$ by
\beql{kappadef}
\kappa_{jk}\colon\ (R_0,\ldots,R_{n-1}) \mapsto (R_0,\ldots,R_{j-1},R_jK_k,R_{j+1},\ldots,R_{n-1}) =: (S_0,\ldots,S_{n-1}).
\eeq
This produces a new group with generators $S_0,\ldots,S_{n-1}$.  We abbreviate $\kappa_{jj}$ to $\kappa_j$, because this is the most important case (and here usually only with $j = 0,1$), but $\kappa_{02}$ is also useful.  Thus $\kappa_j$ interchanges the two possibilities for $R_j$ which can occur in \thmref{fullrankmir}.  Just as with the Petrie operation, though, it must be emphasized that it is by no means generally the case that $\kappa_{jk}$ will yield a C-group when it is applied to another;  for example, for $S_j$ to be an involution, we need $j = k$ or $j \leq k - 2$.  Observe also that $K_{n-1} = R_{n-1}$, so that the Petrie operation of \eqref{petrie} can be written as $\pi = \kappa_{n-3,n-1}$.  

One result, for which we only have a case-by-case (but not a general) proof, is the following.
\bthml{kappa0}
If $P$ is a finite regular polytope of full rank, then $P^{\kappa_0}$ is also a finite regular polytope of full rank.
\ethm

It is instructive to see how the operation $\kappa_0$ acts geometrically on simple examples.  In fact, $\kappa_0$ may do one of three things, even when the original group $\CG$ is a hyperplane reflexion group:  it may double the order, leave it the same, or even halve it.  To illustrate this, in $\BE^3$ take, respectively, the (group of the) tetrahedron, octahedron and cube; note that, in each case, whereas the old facets were of full rank, the new ones (of the polyhedron associated with the new group) are skew polygons, and so are not.  In the planar case, we have $\{p\}^{\kappa_0} = \{q\}$, where $\frac{1}{p} + \frac{1}{q} = \frac{1}{2}$.

\breml{kappamix}
If $K_k \in \scl{R_j,\ldots,R_{n-1}}$, then $\kappa_{jk}$ results in a mixing operation.  
\erem

It would be inappropriate to reproduce all the details of \cite{m_rpfr} here, even in outline form.  However, let us note a few of the salient facts.   We shall say more about three and four dimensions in later sections;  from five dimensions on, things settle in a common pattern.  Recall our conventions that ``regular (or chiral) polytope" will mean ``faithfully realized finite abstract regular (or chiral) polytope'', while ``regular (or chiral) apeirotope" will mean ``discrete faithfully realized abstract regular (or chiral) apeirotope''.

For the regular $n$-polytopes in $\BE^n$, we add to the simplex, cross-polytope and cube the results of applying $\kappa_0$ to each.  From the $n$-simplex $\{3^{n-1}\}$ we obtain a polytope $\{3^{n-1}\}^{\kappa_0}$ with $2(n+1)$ vertices, those of the simplex and its dual;  its group is $S_{n+1} \times C_2$.  For the $n$-cross-polytope, $\{3^{n-2},4\}^{\kappa_0}$ has the same vertices and symmetry group as $\{3^{n-2},4\}$.  With the $n$-cube $\{4,3^{n-2}\}$, there is a distinction between even and odd dimensions $n$.  When $n$ is even, $\{4,3^{n-2}\}^{\kappa_0}$ has the same vertices and symmetry group;  however, when $n$ is odd, $\{4,3^{n-2}\}^{\kappa_0}$ is isomorphic to the \emph{half-cube} $\{4,3^{n-2}\}/2 \cong \{4,3^{n-2}\}_n$, obtained by identifying opposite vertices of the cube.

For the regular $(n+1)$-apeirotopes in $\BE^n$, we can apply the ``apeir'' construction to each of the six $n$-polytopes of the last paragraph.  We also have $\{4,3^{n-2},4\}$, the tiling of space by $n$-cubes, and, finally, $\{4,3^{n-2},4\}^{\kappa_1}$, which is obtained from it by replacing its vertex-figure $\{3^{n-2},4\}$ with $\{3^{n-2},4\}^{\kappa_0}$.  This last is very interesting;  its $3$-face is the Petrie-Coxeter apeirohedron $\{4,6 \hole 4\}$, and, more generally, its facet is the $n$-face of $\{4,3^{m-2},4\}^{\kappa_1}$ for each $m \geq n$.

The following table lists the numbers of regular polytopes and apeirotopes of full rank, according to dimension.

\begin{center}
\begin{tabular}{||c|c|c||}
\hline 
dimension & polytopes & apeirotopes \\
\hline \hline
0 & 1 & - \\
1 & 1 & 1 \\
2 & $\infty$ & 6 \\
3 & 18 & 8 \\
4 & 34 & 18 \\
$\geq 5$ & 6 & 8 \\
\hline 
\end{tabular}
\end{center}

We end the section by quoting another result from \cite{m_rpfr}.  If equality occurs in \propref{rankchirreal}, then (as before) we say that $P$ is \emph{of full rank}.  This result shows that including chiral polytopes does not add any new examples to the previous classification.
\bthml{fullrankchir}
There are no chiral realizations of polytopes of full rank.
\ethm

\section{Regular polytopes in three dimensions}
\label{regpol3}

The paper \cite{ms_rpo} was devoted to the complete classification of the regular polytopes and apeirotopes in $\BE^3$, and so we confine ourselves here to the briefest mention of the techniques employed.

With rank at most $2$, we have the segment in rank $1$, and the polygons (planar and zigzag) and apeirogons (linear, zigzag and helical) in rank $2$.  We say no more about them.

With rank $3$, we first note that the three regular planar tessellations and their Petrials are planar.  There are nine ``classical'' regular polyhedra (the so-called Platonic solids and the Kepler-Poinsot polyhedra -- see \cite[Section~1A]{arp} for discussion of truer attributions), and nine others, which (as a family) can be regarded either as their Petrials, or as the result of applying $\kappa_0$ to them.  There are twelve apeirohedra which are blends of the six planar ones with a segment or apeirogon, and twelve others which are pure (unblended);  of these, except for the Petrie-Coxeter apeirohedra of \cite{c_rsp}, all but one were found by Gr\"unbaum \cite{g_on}, and the last was discovered by Dress \cite{d1,d2}.

The last case of the twelve pure apeirohedra is possibly the most interesting, at least for the methods employed.  A geometric discussion shows that the possible dimension vectors (of the mirrors of the generating reflexions) are given by $(2,1,2),\ (1,1,2),\ (1,2,1)$ and $(1,1,1)$.  If these mirrors are $R_0,R_1,R_2$ (we assume that the initial vertex is $o$, so that $R_1,R_2$ are linear mirrors), define $S_0'$ to be the translate of $R_0$ through $o$, $S_j' := R_j$ for $j = 1,2$, and finally $S_j := S_j'$ or $-S_j'$, according as $R_j$ is a plane or line.  This relates the original symmetry group to one of the crystallographic Coxeter groups $[3,3], [3,4]$ or $[4,3]$ (we need both the latter forms) or the corresponding regular polyhedra;  then the three groups, each with four dimension vectors, result in the twelve apeirohedra.

These apeirohedra are listed in the following table;  for any notation not introduced hitherto, we refer to \cite{ms_rpo} or \cite[Section~7E]{arp}.
\begin{center}
\begin{tabular}{|c||ccc|}
\hline
& $\{3,3\}$ & $\{3,4\}$ & $\{4,3\}$  \\
\hline  \hline
(2,1,2) & $\{6,6 \hole 3\}$ & $\{6,4 \hole 4\}$ & $\{4,6 \hole 4\}$  \\
(1,1,2) & $\{\infty,6\}_{4,4}$ & $\{\infty,4\}_{6,4}$ & $\{\infty,6\}_{6,3}$ 
\\
(1,2,1) & $\{6,6\}_{4}$ & $\{6,4\}_{6}$ & $\{4,6\}_{6}$  \\
(1,1,1) & $\{\infty,3\}^{(a)}$ & $\{\infty,4\}_{\cdot,*3}$ &
$\{\infty,3\}^{(b)}$  \\
\hline
\end{tabular}
\end{center}
The entries in the left column are the dimension vectors $(\dim R_{0},\dim R_{1},\dim R_{2})$, and the remaining columns are indexed by the corresponding finite regular polyhedra.  Of these twelve apeirohedra, nine occur naturally as distinguished members in large families of polyhedra (generally apeirohedra), in which all but two polyhedra are chiral (the two exceptional polyhedra are regular);  we elaborate on this in Section~\ref{chirpol3}.

Finally, there are eight regular $4$-apeirotopes in $\BE^3$ (see \cite[Section~7F]{arp}).  There is the 
regular tiling $\{4,3,4\}$ of space by cubes, the result $\{\{4,6 \hole 4\},\{6,4\}_3\}$ of applying
$\kappa_1$ (or $\pi$) to it, and six more obtained by applying the ``apeir'' operation to the six rational
regular polyhedra, namely, the tetrahedron, octahedron and cube and their Petrials.
\smallskip

\section{Chiral polytopes in three dimensions}
\label{chirpol3}

We now proceed with the enumeration of the (discrete and faithful) chiral polyhedra in $\BE^3$, following \cite{s_cp1,s_cp2}.  Again, we shall not go into details and therefore only briefly summarize the results.

The symmetry group $\CG := \CG(P)$ of a chiral polyhedron $P$ has two orbits on the flags, such that adjacent flags are in distinct orbits. If $\CP$ is the underlying abstract polyhedron, then $\CG$ is isomorphic to $\Ga(\CP)$ or $\Ga^+(\CP)$ according as $\CP$ is chiral or regular. In either case, $\CG = \scl{S_1,S_2}$, where $S_1,S_2$ are the distinguished generators of $\CG$ associated with a base flag $\Ph$ of $P$ and corresponding to the generators $\sg_1,\sg_2$ of $\Ga(\CP)$ or $\Ga^+(\CP)$, respectively.  If $P$ is of type $\{p,q\}$, then
\[ S_1^p = S_2^q = (S_1S_2)^2 = I , \]
but in general there are also other independent relations. If $\Ph$ is replaced by $\Ph^2$ (the adjacent flag with a different $2$-face), then the new pair of generators of $\CG$ are $S_1S_2^2,S_2^{-1}$.  Thus $S_1,S_2$ and $S_1S_2^2,S_2^{-1}$ are the pairs of generators representing the two enantiomorphic forms of $P$. 

As we remarked in Section~\ref{real}, a chiral polyhedron $P$ can be obtained from a variant of Wythoff's construction, applied to a group $\CG=\scl{S_1,S_2}$ with initial vertex a point $v$ fixed by $S_2$ (but not $S_1$).  If we set $T := S_1S_2$, which must be a reflexion in a line or plane, then the base vertex, edge and facet of $P$ are $v$, $v\scl{T}$ and $(v\scl{T})\scl{S_1}$, respectively;  as usual, the other vertices, edges and facets are their images under $\CG$. 

The first step is to determine the possible special groups and their generators.  Recall that, if $R\colon x\mapsto xR' + t$ is a general element of $\CG$, with $R' \in \RO_3$, the orthogonal group, and $t \in \BE^{3}$ a translation vector, then the linear mappings $R'$ form the \emph{special group} $\CG_0$ of $\CG$.  In the present context, $\CG$ must be a crystallographic group in $\BE^3$ and $\CG_0 = \scl{S_1',S_2'}$ a finite subgroup of $\RO_3$.  If $T(\CG)$ denotes the subgroup of all translations in $\CG$, then $\CG_0 \cong \CG/T(\CG)$.  It turns out that the only possible special groups are $[3,3]$ and $[3,4]$ (possibly as $[4,3]$), the full tetrahedral and octahedral group, respectively, and their rotation subgroups $[3,3]^+$ and $[3,4]^+$ (possibly as $[4,3]^+$), as well as the group $[3,3]^*$ obtained from $[3,3]^+$ by adjoining the central inversion in the invariant point of $[3,3]^+$. In particular, this limits the possible Schl\"afli types to $\{4,6\}$, $\{6,4\}$, $\{6,6\}$, $\{\infty,3\}$ and $\{\infty,4\}$.

A chiral polyhedron in $\BE^3$ cannot be finite (by Theorem~\ref{fullrankchir}) or be a blend (its group must be affinely irreducible).  Thus each chiral polyhedron is infinite and pure.  

The possible apeirohedra fall into six infinite $2$-parameter families (up to congruence).  In each family, all but two polyhedra are chiral; the two exceptional polyhedra are regular.  The following table lists the families of polyhedra by the structure of their special group, along with the two regular polyhedra occurring in each family; in three families, one exceptional polyhedron is finite. 
\medskip
\begin{center}
\begin{tabular}{|c|c|c|c|c|c|}
\hline
$[3,3]^*$ & $[4,3]$ & $[3,4]$ & $[3,3]^+$ & $[4,3]^+$ & $[3,4]^+$ \\
\hline\hline
$P(a,b)$ & $Q(c,d)$ & $Q(c,d)^*$ & $P_1(a,b)$ & $P_2(c,d)$ & $P_3(c,d)$ \\
\hline
$\{6,6\}_{4}$ & $\{4,6\}_{6}$ & $\{6,4\}_{6}$ & $\{\infty,3\}^{(a)}$ &
$\{\infty,3\}^{(b)}$ & $\{\infty,4\}_{\cdot,*3}$ \\
$\{6,6 \hole 3\}$ & $\{4,6 \hole 4\}$ & $\{6,4 \hole 4\}$ & $\{3,3\}$ & $\{4,3\}$ &  $\{3,4\}$ \\
\hline
\end{tabular}
\end{center}
\medskip
The columns are indexed by the special groups to which the respective polyhedra correspond;  some groups occur twice but with different pairs of  generators. The second row contains the six families; as we said before, possibly with one exception, all members of a family are apeirohedra. For the first three families, discreteness forces the parameter pairs $a,b$ and $c,d$, respectively, to be relatively prime integers; however, for the last three families, the parameters can be reals. (Thus, when the polyhedra are considered up to similarity, there is a single rational or real parameter, as appropriate.)  

In particular, the chiral polyhedra $P(a,b)$, $Q(c,d)$ and $Q(c,d)^*$ (the dual of $Q(c,d)$) have finite skew faces and skew vertex-figures, and are of types $\{6,6\}$, $\{4,6\}$ or $\{6,4\}$, respectively; remarkably, in each family, the two regular polyhedra have planar faces or vertex-figures.  Recall that no regular polyhedron has finite skew faces and skew vertex-figures (see \cite[Section~7E]{arp}).  On the other hand, the polyhedra $P_1(a,b)$, $P_2(c,d)$ and
$P_3(c,d)$ have infinite faces consisting of helices over triangles, squares or triangles, respectively, and are of types $\{\infty,3\}$, $\{\infty,3\}$ or $\{\infty,4\}$.

The last two rows of the table comprise nine of the twelve pure regular apeirohedra in $\BE^{3}$, namely those listed in the table of Section~\ref{regpol3} with dimension vectors $(1,2,1)$, $(1,1,1)$ or $(2,1,2)$, as well as the three (finite) ``crystallographic'' Platonic polyhedra. The three remaining pure regular apeirohedra $\{\infty,6\}_{4,4}$, $\{\infty,4\}_{6,4}$ and
$\{\infty,6\}_{6,3}$ all have dimension vector $(1,1,2)$ and do not occur in families alongside chiral polyhedra. 

We now display the families of polyhedra, with the various known relationships among them.  These complement the known relationships between regular polyhedra (see \cite[Section~7E]{arp}).  Three operations on (chiral or regular) polyhedra and their groups $\CG$ are involved:\  the duality operation $\delta$, the second facetting operation $\vp_2$ and the halving operation $\eta$ (see Section~\ref{regchirpol} or \cite{s_cp1}). In terms of the generators of $\CG$ they are defined as follows:
\[\begin{array}{rlll}
\delta\colon   & (S_1,S_2)  &\mapsto& (S_2^{-1},S_1^{-1}) ,\\
\vp_2\colon & (S_1,S_2)   &\mapsto& (S_1S_2^{-1},S_2^2), \\
\eta\colon      & (S_1,S_2) &\mapsto& (S_1^2S_2,S_2^{-1}). 
\end{array} \]
In each case, the pair of elements on the right are the generators for the group of a new polyhedron, namely the image of the given polyhedron under $\delta$, $\vp_2$ or $\eta$, respectively.  

The following diagram emphasizes operations relating families rather than individual polyhedra. In particular, we drop the parameters from the notation;  for example, $P_1$ denotes the family of polyhedra $P_1(a,b)$.  
\beq
\label{displayone}
\begin{matrix}
Q^* & \stackrel{\delta}{\longleftrightarrow} \mkern-30mu & Q & \mkern-30mu 
\stackrel{\vp_2}{\longrightarrow}& P_2 \cr \cr
& &\;\;\downarrow\! {\scriptstyle\eta} & & & \cr \cr
P_3 & & P & \mkern-30mu \stackrel{\vp_2}{\longrightarrow}& P_1 \cr
& & 
\begin{picture}(60,60)
\put(9,46){\oval(42,42)[b]}
\put(9,46){\oval(42,42)[tl]}
\put(9,67){\vector(1,0){2}}
\put(-19,34){$\scriptstyle\delta$}
\end{picture}
&&
\end{matrix} 
\eeq
Observe that, in the diagram, $P_3$ is not connected to any other family;  it is an interesting open question if there exists a relationship between $P_3$ and any other family.  The circular arrow in the diagram indicates the self-duality of the family (in fact, of each of its polyhedra). The operations $\delta$ and $\vp_2$ map a polyhedron to one with the same parameter pair, either $a,b$ or $c,d$. However, $\eta$ replaces $c,d$ by the new pair $c-d,c+d$.  Moreover, note that $\vp_2$, when applicable, maps a polyhedron to one in the same row of the table.

For a discussion of other classes of highly symmetric polyhedra in $\BE^3$ we refer the reader to, for example, \cite{g_ap}.

\section{Regular polytopes in four dimensions}
\label{regpol4}

Just as is the case with the classical regular polytopes and apeirotopes, the richest family of full rank occurs in $\BE^4$.  Again, we do not wish to go into the results of \cite{m_rpfr} in great detail;  instead, we shall concentrate on a few plums.

We have already accounted for the effects of $\kappa_0$;  we merely note that the sixteen classical regular (convex and star-) $4$-polytopes give rise to another sixteen in this way.  However, we can also apply $\pi$ to the $4$-cube $\{4,3,3\}$, to obtain
\[  \{4,3,3\}^\pi = \{\{4,4 \hole 4\},\{4,3\}_3\}.  \]
That is, the facets are toroids, and the vertex-figure is the half-$3$-cube;  moreover, the polytope is universal of this kind.  The final instance (of the $34$ in the table of \sectref{fullrank}) is obtained by applying $\kappa_0$ to this last.

These two finite polytopes just mentioned contribute two regular $5$-apeirotopes via the ``apeir" construction.  Leaving aside the examples already discussed in \sectref{fullrank}, then for the $5$-apeirotopes there remain those obtained from $\{3,3,4,3\}$ and its dual $\{3,4,3,3\}$.  For the first, we can apply $\kappa_1$ (that is, apply $\kappa_0$ to its vertex-figure $\{3,4,3\}$);  we get an apeirotope whose $3$-faces are Petrie-Coxeter apeirohedra $\{6,6 \hole3\}$.  For the other, we can first apply $\kappa_1$;  the resulting apeirotope has $3$-faces the last Petrie-Coxeter apeirohedron $\{6,4 \hole 4\}$.  To both of these (that is, $\{3,4,3,3\}$ and $\{3,4,3,3\}^{\kappa_{1}}$), we can now apply $\pi$ as well;  the $3$-faces remain as they were (that is, octahedra $\{3,4\}$ or $\{6,4 \hole 4\}$, respectively);  the facet of the first is the universal apeirotope $\{\{3,4\},\{4,4 \hole 4\}\}$ (we comment on this further in \sectref{opprob}).

It is a striking fact that all three Petrie-Coxeter apeirohedra in $\BE^3$ occur as $3$-faces of regular $5$-apeirotopes in $\BE^4$ (one of them twice).

We next discuss the recent (as yet unpublished) classification of the four-dimensional (finite) regular polyhedra.  Those polyhedra with planar faces were all found in \cite{abm,bracho};  the methods we employ in \cite{m_fdrp} are akin to those used in \cite{m_rpfr,ms_rpo}, and are, we feel, much simpler.

As we have already pointed out in \sectref{real}, our strategy is to determine what possible dimension vectors can occur, and then to enumerate every polytope in the corresponding subclasses.  \thmref{dimmirror} provides a starting point;  in the current case, the dimension vector must satisfy
\[  \dim R_0 \geq 1, \quad \dim R_1 \geq 2, \quad \dim R_2 \geq 2.  \]

We now proceed as follows.  As essentially the same trick we perform in $\BE^3$, if the mirror $R_0$ satisfies $\dim R_0 = 1$, then we can replace it by
\[  -R_0 = R_0^\perp,  \]
its orthogonal complement, which (as an isometry) is its product with the central inversion $-I$;  we refer to this more general operation as $\kappa_0$ as well.  We always obtain another finite group $\CG'$;  in fact,
\[  |\CG'| = \half|\CG|,\ |\CG|,\ \mbox{or } 2|\CG|.  \]

Next, if $\dim R_0 = 2$ and $\dim R_2 = 3$ (or vice versa, but this case will have to be excluded), then we can replace $R_0$ by $R_0R_2$, that is, apply (or reverse) the Petrie operation $\pi$;  bearing in mind \eqref{prodaxis}, the new $R_0$ has $\dim R_0 = 1$ or $3$, and in the former case we can proceed as previously.  

Finally, as long as our (possibly new) group contains a hyperplane reflexion (that is, $\dim R_j = 3$ for some $j$), we can regard $\CG$ as a reflexion (Coxeter) group, on which certain involutions with $2$-dimensional mirrors act as automorphisms (more precisely, $\CG$ is the corresponding semi-direct product).  When we have carried out the foregoing procedures, only the dimension vectors $(3,2,3)$ and $(2,3,2)$ need to be considered.  For classification purposes, we then reverse the procedure:  the starting point is a Coxeter group, not necessarily with standard
generators, which can be represented by a diagram that permits permutation of its nodes.

We give a couple of simple examples of what happens in the cases $(3,2,3)$ and $(2,3,2)$ in a little detail, and then comment on the remaining cases (with the exception of $(2,2,2)$) more briefly.  We list them according to their dimension vectors.

\bitem

\item $(3,2,3)$:  from the group $[3,4,3]$ of the regular $24$-cell, we derive the diagrams
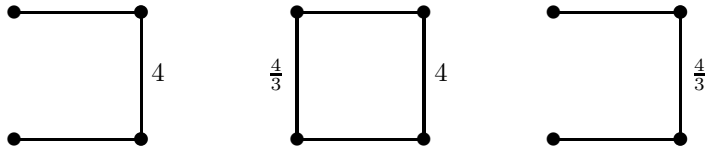
\begin{center}
\bpc(250,70)
\Horone{0,11}
\Horone{0,59}
\Verone{48,11}
\put(52,33){$4$}
\Resq{107,11}
\put(96,33){$\frac{4}{3}$}
\put(159,33){$4$}
\Horone{204,11}
\Horone{204,59}
\Verone{252,11}
\put(256,33){$\frac{4}{3}$}
\epc
\end{center}
each of which permits a top-to-bottom flip, and thereby gives two dual regular polyhedra with dimension vectors $(3,2,3)$.  (From the first diagram, we obtain the polyhedra $\{4,8 \hole 3\}$ and $\{8,4 \hole 3\}$ of \cite{c_rsp}.)  Similar examples derive from the diagram
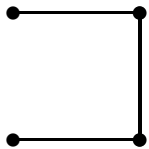
\begin{center}
\bpc(48,70)
\Horone{0,11}
\Horone{0,59}
\Verone{48,11}
\epc
\end{center}

\item $(2,3,2)$:  the general case is derived from a diagram 
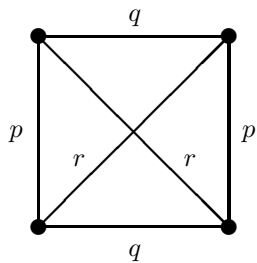
\begin{center}
\bpc(110,100)
\REsq{19,14}
\thicklines
\put(19,14){\line(1,1){72}}
\put(19,86){\line(1,-1){72}}
\multiput(8,48)(88,0){2}{$p$}
\multiput(32,37)(42,0){2}{$r$}
\multiput(53,3)(0,89){2}{$q$}
\epc
\end{center}
with horizontal and vertical flips.  This gives rise to a polyhedron of type $\{2p,2q\}_{2r}$, from which are obtained up to five others by duality and Petriality.  (There is a restriction on $q$:  it must not be a fraction with even denominator.)  As a specific instance, the full family of six is obtained when $\{p,q,r\} = \{3,4,\frac{4}{3}\}$.

\item $(3,3,3)$:  this corresponds to three-dimensional polyhedra, and so is excluded (but only on these grounds).  

\item $(1,3,3)$:  this is allowed;  $\kappa_0$ can be applied to the case $(3,3,3)$.  

\item $(2,3,3)$:  this is obtained from $(3,3,3)$ or $(1,3,3)$ by Petriality;  therefore, the first possibility must be excluded. 

\item $(3,3,2)$:  this would be obtained from $(2,3,3)$ by duality;  however, in the allowed case, the faces of the original are centred at $o$, and so the dual must be excluded.

\item $(1,3,2)$:  this would be obtained from $(3,3,2)$ by applying $\kappa_0$, and so is also disallowed. 

\item $(1,2,3)$:  this arises from $(3,2,3)$ by applying $\kappa_0$.

\item $(2,2,3)$:  this is obtained from $(3,2,3)$ or $(1,2,3)$ by Petriality.  

\item $(3,2,2)$:  this would arise from $(2,2,3)$ by duality.  However, it may be seen that (with either possibility) the product $R_0R_1$ of the corresponding reflexions $R_0$ and $R_1$ in the original is a double rotation (in two orthogonal planes), since $R_0 \cap R_1 = \{o\}$;  it follows that the class cannot occur.  

\item $(1,2,2)$:  this would be obtained from $(3,2,2)$ by applying $\kappa_0$, and so it too must be excluded. 

\eitem

It is notable that only the groups $[3,3,3]$ and $[3,4,3]$ give rise to polyhedra in the classes $(3,2,3)$ and $(2,3,2)$ and those derived from them.  Even though other finite reflexion groups in $\BE^4$ permit diagram automorphisms (for suitably chosen generators), these are inner, and then the corresponding ``polyhedra'' degenerate.

The anomalous case is dimension vector $(2,2,2)$, to which the notion of a Coxeter group with outer automorphisms is inapplicable.  Indeed, some examples of this kind cannot be related to Coxeter groups in any meaningful way.  The approach here is through quaternions.  Each isometry which occurs in such a group is a rotation (that is, lies in $\SO_4$), and so can be represented by a quaternionic transformation of the form
\beql{quatmap}
x \mapsto \ol{a}xb,
\eeq
where $a,\ b$ are unit quaternions (recall that $a^{-1} = \ol{a}$).  In keeping with our usual conventions, mappings are thought of as acting on the right;  thus it must be the inverse of a quaternion which provides an appropriate mapping when acting on the left.  For the mapping \eqref{quatmap} to be a reflexion, both $a$ and $b$ must be pure imaginary.  Our symmetry group $\CG$ gives rise to two groups $\CG_L$ and $\CG_R$ of the left-acting quaternions $a$ and right-acting quaternions $b$;  then $\CG$ is a certain quotient of $\CG_L \times \CG_R$ (for further details at this stage, we refer the reader to \cite{dv}).  Further, there are then quotients $G_L,\ G_R$ of $\CG_L,\ \CG_R$ in $\SO_3$, each by normal subgroups of index $2$, and these
are generated by half-turns about lines in $\BE^3$.  If $a = \cos\vt + u\sin\vt$, with $u$ pure imaginary, then the image of $a$ under the homomorphism from $\CG_L$ to $G_L$ is a rotation through $2\vt$ about the axis in $\BE^3$ through $u$, when the latter is regarded as a unit vector in $\BE^3$.  Thus, for example, if $a$ is pure imaginary, then its image is the half-turn about the axis in $\BE^3$ through $a$;  it is important to note that this half-turn lifts to two pure imaginary quaternions $\pm a$.  The only groups which can occur as such groups $G_L$ or $G_R$ are dihedral, octahedral or icosahedral;  the cyclic and tetrahedral groups do not contain enough half-turns.  Finally, if the generating reflexions are 
\[  xR_j := \ol{a}_jxb_j = -a_jxb_j \quad (j = 0,1,2),  \] 
then (as scalar products of vectors in $\BE^3$), 
\[  \scl{a_1,a_2} = \pm\scl{b_1,b_2},  \]
because the product $R_1R_2$ must have a $2$-dimensional axis.  However, the opposite must be true for the product $R_0R_1$, because this has to be a double rotation.

In summary, the following ingredients go into the enumeration.  First, two groups in $\BE^3$ generated by half-turns:  these are a dihedral group $D_{2k}$ ($k$ can only take the values $2$, $3$ or $5$), the octahedral group $S_4 = [3,3] = [3,4]_3$ or the isosahedral group $A_5 = [3,5]_5$.  Second, for $j = 1,2$, two regular polyhedra of type $\{r_j,q\}$ (with the same $q$);  here, we must allow $r_j > 1$, rather than the usual $r_j \geq 2$, to account for two possible liftings of the half-turns contributing to $R_0$.  We then obtain a polyhedron of type $\{p,q\}$, where the face $\{p\}$ is of the form $\{p_1\} \# \{p_2\}$, with
\[  \frac{1}{p_j} = \frac{1}{2}\left( \pm \frac{1}{r_1} \pm \frac{1}{r_2} \right),  \]
where the signs are chosen so that $p_j > 2$ for $j = 1,2$.  It is convenient to write the face, instead, as
\[  \left\{\frac{p}{d_1,d_2}\right\}, \qquad \mbox{with} \quad p_j = \frac{p}{d_j}  \]
(in lowest terms) for $j = 1,2$. 

As a specific example, if $r_1 = 3,\ r_2 = \fh$ and $q = 5$, then we obtain a polyhedron of type
\[  \{\tfrac{30}{1,11},5\}.  \]
However, if we replace $\fh$ by $\frac{5}{3}$ (or $3$ by $\frac{3}{2}$), indicating a different choice of lifting for $R_0$, then we obtain type
\[  \{\tfrac{15}{2,7},5\}.  \]

\brem
A further comment is in order here.  An opposite orthogonal transformation of $\BE^4$ is of the form
\[  x \mapsto\ol{a}\,\ol{x}b,  \]
with $a,\ b$ as before.  In a group $\CG$ containing such transformations, the corresponding left and right groups $\CG_L$ and $\CG_R$ must be conjugate in the whole group of unit quaternions.  Thus one could also use quaternions to investigate the classes other than $(2,2,2)$;  however, the methods which we have already described are more efficacious.
\erem

\section{Open problems}
\label{opprob}

As the dimension increases, so there are more possibilities for the ranks of faithfully realized regular or chiral polytopes or apeirotopes.  In full rank, the regular cases are classified, and chirality does not occur.  In $\BE^4$, therefore, the open cases are the (finite) chiral polytopes of rank $3$, and the regular or chiral apeirotopes of ranks $3$ and $4$. 

We look at the regular cases first;  we begin with rank $4$.  Each of the eight regular $4$-apeirotopes in $\BE^3$ can be blended (mixed) with a segment or an apeirogon;  this gives $16$ blended examples.  Next, the ``apeir'' construction described at the end of \sectref{real} can be applied to any of the four-dimensional rational regular polyhedra.  Finally, certain of the facets of the regular apeirotopes of full rank in $\BE^4$ are $4$-apeirotopes.  It is possible that there are not too many more examples which do not fall under one of these three categories,
and maybe even none at all.

Incidentally, there is only one four-dimensional $4$-apeirotope whose facets are finite regular polyhedra.  This is the universal $\{\{3,4\},\{4,4 \hole 4\}\}$, with facet the octahedron $\{3,4\}$ and vertex-figure the toroid $\{4,4 \hole 4\}$, which, as noted in \sectref{regpol4},  is the facet of the $5$-apeirotope $\{3,4,3,3\}^\pi$.  (Compare \cite[Theorem~10B3]{arp} with $s = 4$ in the dual form, and the preceding discussion.)  To see that this is the only example, observe that there are no four-dimensional (finite) regular polyhedra with triangular faces (nor
with pentagons or pentagrams either, but these must be excluded on crystallographic grounds).  Hence, the only possible vertex-figure has square faces, which means that the facet must be an octahedron or its Petrial $\{6,4\}_3$.  In turn, the vertex-figure must be a regular polyhedron with square faces, and circumradius equal to its edge-length;  this forces it to be $\{4,4 \hole 4\}$.  Finally, direct calculation shows that, in fact, $\{6,4\}_3$ cannot actually be a facet in such a way.

As for the four-dimensional regular apeirohedra, a mere glance at some of the possibilities shows that the enumeration problem is likely to be rather hard.  For example, in $\BE^2$ the apeirohedron $\{\fh,10\}$ is non-discrete;  however, when it is blended with its isomorphic copy $\{5,\tth\}$ in $\BE^4$, a discrete regular apeirohedron of type $\{5,10\}$ is obtained. Several similar examples also occur.

There are also examples derived from complex reflexion groups in $\BC^2$, which we regard as real groups in $\BE^4$ generated by reflexions with $2$-dimensional mirrors.  A curiosity is the following. We can twist the first of the two diagrams below by the dihedral group $D_3$ (or symmetric group $S_3$), and the second by $C_2$.  We then actually obtain the same geometric group;  however, the outer automorphisms of one correspond to the generating reflexions of the other (and vice versa).  We refer to \cite[Section~9D]{arp} for the background here.
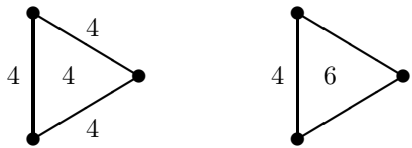
\begin{center}
\bpc(150,68)
\Rtri{10,10}  \Rtri{110,10}
\multiput(30,11)(0,38){2}{$4$}
\put(21,31){$4$}
\put(120,31){$6$}
\multiput(0,31)(100,0){2}{$4$}
\epc
\end{center}

We now turn to chiral polytopes and apeirotopes.  For the latter, various infinite families of chiral apeirohedra were described in \cite{s_cp1,s_cp2} (see Section~\ref{chirpol3});  each such apeirohedron can be blended with a segment or an apeirogon to give a four-dimensional chiral apeirohedron.  Finally, there are plenty of finite chiral polyhedra in $\BE^4$;  for example, each chiral toroid $\{4,4\}_{(s,t)}$ is realizable.  Whether there exist non-toroidal finite chiral polyhedra in $\BE^4$ is a nice open question.

Finally, presentations for the symmetry groups have only been fully worked out for the $3$-dimensional regular polyhedra and apeirotopes (see \cite[Sections~7E, 7F]{arp}). For higher dimensions, presentations are known for certain classes of polytopes, for example, the regular star-polytopes (see \cite[Section 7D]{arp} or \cite{m_grsp}). In this context, the main tool is the so-called ``circuit criterion", which states that the automorphism group of an abstract polytope $\CP$ (and thus the symmetry group of a faithful realization) is determined by the group of its vertex-figure and the circuit structure of the edge-graph of $\CP$ (see \cite[Section 2F]{arp} for more details). A variant of this method should also succeed in the chiral case. In particular, there is an interest in presentations for the symmetry groups of the $3$-dimensional chiral apeirohedra. Here we do not know if the corresponding  abstract apeirohedra are also chiral or if they are regular. Settling this question may have to be the first step in arriving at presentations for their symmetry groups.

\end{document}